\documentclass{amsart}
\usepackage{amssymb}
\usepackage{amsfonts}
\usepackage{amsmath}
\usepackage{amsthm}

\newtheorem{teo}{Theorem}

\theoremstyle{definition}
\newtheorem{definicion}[teo]{Definition}
\newcommand{\cuerpo}[1]{\mathbb{#1}}
\newcommand{\p}{\cuerpo{P}}
\newcommand{\C}{\cuerpo{C}}
\DeclareMathOperator{\rg}{rk}
\title{A characterization of the rational normal curve}

\author{Gonzalo Comas} 
\email{gcomas@dm.uba.ar}
\address{Departamento de Matematica FCEyN-UBA\\
Ciudad Universitaria (1428) Capital Federal\\
Argentina}
\begin{document}

\begin{abstract}
We give a characterization of the rational normal curve in terms of the
rank function associated to a curve.
\end{abstract}

\maketitle

\begin{section}{Introduction}
Let $C \subset \p^n$ be a nondegenerate and nonsingular curve
over $\C$. Let $x \in \p^n$
be a point. We define the $C$-rank of $x$ as the smallest number $r$
such that $x$ lies in the linear span of $r$ points of $C$. In \cite{mio}
we describe the strata of points having constant rank in the case
in which $C$
is a rational normal curve, that is, a degree $n$ curve. It is shown there
that all points belonging to a tangent line to $C$ have rank $n$, excepting
the point of tangency, which has rank one. 
In \cite{comas}, it is shown that if $C$ has positive genus, and the
inmersion is given by a complete linear system of degree $d$ divisors
(with $d \geq 10g-1$),
then all points in $\p^n$ have rank less or equal than $n-g$, and all points
in tangent lines have rank $n-g$ (excepting the points of tangency).

In this article we show that the only curve $C \subset \p^n$ 
such that all points belonging to a tangent line to $C$ have rank $n$ (excepting
the point of tangency), is the rational normal curve. 

{\bf Acknowledgments.} We
would like to thank J.~M.~Landsberg for posing this question.

\end{section}

\begin{section}{Main result.}
\begin{definicion}
Let $C \subset \p^n$ be a nondegenerate and nonsingular curve. Let $x \in \p^n$
be a point. We define the $C$-rank of $x$ (and note it $\rg x$)
as the smallest number $r$
such that $x$ lies in the linear span of $r$ points of $C$.
\end{definicion}

We will prove the following theorem

\begin{teo}\label{principal}
The rational normal curve is the only nondegenerate and nonsingular curve
in $\p^n$ such that all the points in its tangent lines (excepting the point of tangency)
have rank $n$.
\end{teo}

First we show that the rational normal curve has the desired property.

\begin{teo} 
Let $C \subset \p^n$ be the rational normal curve, let $p \in C$
and let $x \in T_p(C)$ a point such that $x \neq p$. Then 
$\rg(x)=n$.
\end{teo}
\begin{proof}
Assume that $r=\rg(x) <n$.
Then $x \in \langle p_1,\dots,p_r\rangle$ for
$p_1,\dots,p_r \in C$.
Let $L$ be the tangent line to $C$ at $p$,
$M$ be the linear span of $p_1,\dots,p_r$
and $N$ be the linear span of $L$ and $M$.
We have that $x \in L \cap M$, therefore
$N$ has dimension
less than or equal to $r$.

If $L\cap M = \{x\}$, that is, if $p \neq p_i \forall i$,
then $\dim N=r$.

If, on the other hand, $p=p_i$ for some $i$, then $L \cap M$
contains the line spanned by $p$ and $x$, and therefore $L \cap M=L$.
Then $\dim N=r-1$.
 
But the intersection of a linear variety
of dimension $k$ cuts the rational normal curve
in a divisor of degree less than or equal to $k+1$. In both cases
the intersection of $N$ with $C$ is a divisor of degree $\dim N +2$.

Therefore $\rg (x) \geq n$.

Now we show that every point in $\p^n$ has rank less than or equal to $n$.
Let us consider the linear system ${\mathcal D}$ of hyperplanes that contain $x$.
It is clear that the intersection of all the members in ${\mathcal D}$
is the set which has $x$ as its only element.
We want to show that one of the members in ${\mathcal D}$ has no multiple points.
If every member has multiple points, then by Bertini's Theorem
there exists a point $p \in C$ that is a multiple point of every
member in ${\mathcal D}$. But then the intersection
of all the hyperplanes contains the tangent line to $C$ at $p$,
which is a contradiction. Therefore
one of the members in ${\mathcal D}$ has no multiple points, and 
the rank of $x$ is less than or equal to $n$.
\end{proof}

\begin{proof}[Proof of Theorem \ref{principal}.]
We already proved that the rational normal curve
has the desired property.
Now we will show that if $C$ is not the rational normal curve,
then there are points on its tangent lines which have neither rank 1 nor
rank $n$.
Since $C$ is not the rational normal curve, the degree of $C$ is greater than $n$.

If $n=2$, a generic tangent line to $C$ cuts $C$ in another point.
Therefore the generic tangent line to $C$ contains a point that is not
the point of tangency, which has rank different from 2.

From now on $n \geq 3$.
Let $p \in C$ be a point.
If the tangent line to $C$ at $p$ cuts $C$ in another point $q$
we are done, since $q$ has rank one.
Since there are only finite points
$p$ such the tangent line to $C$ at $p$  cuts $C$ with multiplicty
greater than 2, we assume that the intersection
of the tangent line to $C$ at $p$ and $C$ is the divisor $ 2p $.

Let us consider first the case $n\geq 4$. Let $L$ be the tangent line 
to $C$ at $p$ and let $\pi_L: \p^n \setminus L \to \p^{n-2}$
be the projection with center $L$. 
The restriction of $\pi_L$ to $C\setminus \{p\}$
extends to a function on all $C$. Let $C'$
be the image $\pi_L(C)$, which is a
nondegenerated degree $d-2$ curve ($d-2 \geq 3$). By the general position
theorem (\cite{acgh}, pag.~109), the generic hyperplane section to $C'$ cuts $C'$
in $d-2$ disctint points such that any $n-2$ are in general position.
Furthermore, we can assume that the generic hyperplane section does
not contain any of the singular points of $C'$. We can also assume
that the generic hyperplane section does not contain the image of $p$
by $\pi_L$. 
Let $q_1,\dots,q_{d-2} \in C'$ be the points of intersection of $C'$
with a general hyperplane. 
By the assumptions we made 
we can choose $p_1,\dots,p_{d-2} \in C$ (with $p_i \neq p_j$)
such that $\pi_L(p_i)=q_i$.

We now use the fact that there is
a bijection between hyperplanes section of $C'$ and hyperplanes
sections to $C$ that contain $L$. 
By this bijection the hyperplane  $\langle q_1 , \dots, q_{d-2}\rangle$
corresponds to the hyperplane $\langle L,p_1,\dots,p_{d-2}\rangle$.
Let $\Lambda \subset \p^n$ be the linear variety 
$\Lambda = \langle p_1,\dots,p_{d-2}\rangle$. As $\pi_L(\Lambda)$
is a variety of dimension $n-3$, the possible values of $\dim \Lambda$
are $n-3,n-2$ or $n-1$.
If $\dim \Lambda = n-2$, then the intersection $L \cap \Lambda$
is a point $x \not \in C$. We can choose $n-1$ of the $p_i$'s such that 
they generate $\Lambda$. Then we have $\rg x \leq n-1$.
If $\dim \Lambda =n-1$ then $L \subset \Lambda$. Let us choose $n$ of 
the $p_i$'s that generate $\Lambda$. The linear variety
generated by any $n-1$ of these points cuts $L$ in a point $x$ not in $C$. 
Then $\rg x \leq n-1$.

To finish the demonstration we have to prove that for the
generic hyperplane
section $H'$ of $C'$ the corresponding $\Lambda$ has dimension $n-1$
or $n-2$.
So let us assume that $\dim \Lambda=n-3$. As $d-2 > n-2$,
we have that $\Lambda$ is a $(n-3)$-plane that cuts
$C$ in a divisor of degree greater than $n-2$.
Now we compare dimensions.
There are $\infty^{n-2}$
generic hyperplane sections of $C'$.
On the other hand, a nonsingular and nondegenerate curve
can have at most $\infty^{n-3}$ $(n-1)$-secant $(n-3)$-planes
(\cite{acgh}, pag.~152).
Therefore, for the generic hyperplane section $H'$ 
we must have $\dim \Lambda=n-1$ or $n-2$.

If $n=3$ we study the family of planes containing $L$.
If a plane $H$ cuts $C$ in less than $d-2$ points,
then there are double points in this intersection. That means
that $H$ contains another tangent line $L'$
that must cut $L$. Since two generic tangent lines to $C$
do no meet, for general $H$ containing
$L$ the intersection of $H$ and $C$ is the divisor $2p + p_1+\dots+p_{d-2}$
with $p_i \neq p \; \forall i$.
We consider the linear variety $\Lambda=\langle p_1 + \dots + p_{d-2}\rangle$.
If $\dim \Lambda=1$, then the intersection
of $L$ and $\Lambda$ is a point that must have rank $2$
since $L$ does not cut $C$ away from $p$.
If on the other hand $\dim \Lambda=2$,
we can choose $3$ of the $p_i$'s
that generate $\Lambda$. The intersection
of $L$ with one of the lines generated by two of the chosen $p_i$'s
is a 2 rank point.
\end{proof}

\end{section}

%\nocite{hart}
\nocite{acgh}
\nocite{comas}
\bibliographystyle{amsalpha} 
\bibliography{biblio} 

\providecommand{\bysame}{\leavevmode\hbox to3em{\hrulefill}\thinspace}
\providecommand{\MR}{\relax\ifhmode\unskip\space\fi MR }
% \MRhref is called by the amsart/book/proc definition of \MR.
\providecommand{\MRhref}[2]{%
  \href{http://www.ams.org/mathscinet-getitem?mr=#1}{#2}
}
\providecommand{\href}[2]{#2}
\begin{thebibliography}{ACGH85}

\bibitem[ACGH85]{acgh}
E.~Arbarello, M.~Cornalba, P.~A. Griffiths, and J.~Harris, \emph{Geometry of
  algebraic curves}, vol.~1, Springer Verlag, 1985.

\bibitem[Com07]{comas}
Gonzalo Comas, \emph{The rank associated to a projective curve}, Preprint,
  2007.

\bibitem[CS02]{mio}
Gonzalo Comas and Malena Seiguer, \emph{On the rank of a binary form},
  Preprint, 2002.

\end{thebibliography}

\end{document}